\documentclass{amsart}
\usepackage{amsmath,amssymb}
\usepackage[dvips]{graphics}
\usepackage[all]{xy}

\newtheorem{Theorem}{Theorem}[section]

\newtheorem{Proposition}[Theorem]{Proposition}

\newtheorem{Example}[Theorem]{Example}

\newtheorem{Remark}[Theorem]{Remark}

\makeatletter
\@addtoreset{figure}{section}
\def\@thmcountersep{-}
\makeatother

\numberwithin{equation}{section}

\begin{document}

\title{Discontinuous maps whose iterations are continuous}

\author{Kouki Taniyama}
\address{Department of Mathematics, School of Education, Waseda University, Nishi-Waseda 1-6-1, Shinjuku-ku, Tokyo, 169-8050, Japan}
\email{taniyama@waseda.jp}
\thanks{The author was partially supported by Grant-in-Aid for Scientific Research (C) (No. 24540100), Japan Society for the Promotion of Science.}


\date{}

\dedicatory{}


\begin{abstract}
Let $X$ be a topological space and $f:X\to X$ a bijection. Let ${\mathcal C}(X,f)$ be a set of integers such that an integer $n$ is an element of ${\mathcal C}(X,f)$ if and only if the bijection $f^n:X\to X$ is continuous. A subset $S$ of the set of integers ${\mathbb Z}$ is said to be realizable if there is a topological space $X$ and a bijection $f:X\to X$ such that $S={\mathcal C}(X,f)$. A subset $S$ of ${\mathbb Z}$ containing $0$ is called a submonoid of ${\mathbb Z}$ if the sum of any two elements of $S$ is also an element of $S$. We show that a subset $S$ of ${\mathbb Z}$ is realizable if and only if $S$ is a submonoid of ${\mathbb Z}$. Then we generalize this result to any submonoid in any group. 
\end{abstract}

\maketitle

\section{Introduction} 

Let $X$ be a topological space and $f:X\to X$ a bijection. By $f^{-1}:X\to X$ we denoted the inverse mapping of $f$. For each integer $n$ we define a bijection $f^n:X\to X$ by
\[
f^n=\left\{ \begin{array}{ll}
\underbrace{f\circ f\circ\cdots\circ f}_n & (n>0)\\ \\
{\rm id}_X & (n=0)\\ \\
\underbrace{f^{-1}\circ f^{-1}\circ\cdots\circ f^{-1}}_{-n} & (n<0).\\
\end{array} \right.
\]
We note that $f^n\circ f^m=f^{m+n}$ for any integers $m$ and $n$. 
Let ${\mathbb Z}$ be the set of all integers. We define a subset ${\mathcal C}(X,f)$ of ${\mathbb Z}$ by
\[
{\mathcal C}(X,f)=\{n\in{\mathbb Z}|f^n:X\to X{\rm \ is\ continuous.}\}.
\]
A subset $S$ of ${\mathbb Z}$ is said to be {\it realizable} if there is a topological space $X$ and a bijection $f:X\to X$ such that $S={\mathcal C}(X,f)$. A subset $S$ of ${\mathbb Z}$ is called a {\it submonoid} of ${\mathbb Z}$ if $S$ satisfies the following two conditions. 

\begin{enumerate}
\item[(1)]$S$ contains $0$, 
\item[(2)]if $S$ contains $a$ and $b$ then $S$ contains $a+b$. 
\end{enumerate}

\noindent
Note that it is not necessary that $S$ contains $a-b$. 

\vskip 3mm

\begin{Example}\label{submonoid}
{\rm The following subsets of ${\mathbb Z}$ are submonoids of ${\mathbb Z}$. 

\noindent
${\mathbb Z}$, $\{n\in{\mathbb Z}|n\geq0\}$, $\{0\}\cup\{n\in{\mathbb Z}|n\leq-3\}$, $\{2n|n\in{\mathbb Z}\}$, $\{0\}\cup\{3n|n\in{\mathbb Z}, n\geq2\}$, $\{3a+5b|a,b\in{\mathbb Z}, a,b\geq0\}=\{0,3,5,6,8\}\cup\{n\in{\mathbb Z}|n\geq9\}$, $\{0\}$. 

}
\end{Example}

\vskip 3mm

\begin{Theorem}\label{main-theorem}
A subset $S$ of the set of all integers ${\mathbb Z}$ is realizable if and only if $S$ is a submonoid of ${\mathbb Z}$. 
\end{Theorem}

\vskip 3mm

We generalize Theorem \ref{main-theorem} to any submonoid in any group in the third section. 

\section{Proof of Theorem \ref{main-theorem}} 

\begin{Proposition}\label{proposition}
Let $X$ be a topological space and $f:X\to X$ a bijection. Then the subset ${\mathcal C}(X,f)$ of ${\mathbb Z}$ is a submonoid of ${\mathbb Z}$. 
\end{Proposition}

\vskip 3mm

\noindent{\bf Proof.}
Since $f^0={\rm id}_X$ is continuous the set ${\mathcal C}(X,f)$ contains $0$. Suppose that ${\mathcal C}(X,f)$ contains $a$ and $b$. Then $f^a$ and $f^b$ are continuous. Then the composition $f^b\circ f^a=f^{a+b}$ is also continuous. Therefore ${\mathcal C}(X,f)$ contains $a+b$. $\Box$

\vskip 3mm

\noindent{\bf Proof of Theorem \ref{main-theorem}.}
It follows from Proposition \ref{proposition} that if $S$ is realizable then $S$ is a submonoid of ${\mathbb Z}$. We will show that if $S$ is a submonoid of ${\mathbb Z}$ then $S$ is realizable. 
Let $S$ be a submonoid of ${\mathbb Z}$. 
For each integer $n$ we define a subset $X_n$ of the $2$-dimensional Euclidean space ${\mathbb R}^2$ as follows. 
\[
X_n=\left\{ \begin{array}{ll}
\{n\}\times[0,2) & (n\in S)\\
\{n\}\times([0,1)\cup[2,3)) & (n\in ({\mathbb Z}\setminus S)).\\
\end{array} \right.
\]
Let $\displaystyle{X=\bigcup_{n\in{\mathbb Z}}X_n}$. Then $X$ is a topological subspace of ${\mathbb R}^2$. 
Let $f:X\to X$ be a bijection defined by the followings. 

\begin{enumerate}
\item[(1)]if $n,n+1\in S$, then $f((n,x))=(n+1,x)$ for each $x\in[0,2)$,
\item[(2)]if $n,n+1\in({\mathbb Z}\setminus S)$, then $f((n,x))=(n+1,x)$ for each $x\in([0,1)\cup[2,3))$, 
\item[(3)]if $n\in S$ and $n+1\in({\mathbb Z}\setminus S)$, then $f((n,x))=(n+1,x)$ for each $x\in[0,1)$ and $f((n,x))=(n+1,x+1)$ for each $x\in[1,2)$, 
\item[(4)]if $n\in({\mathbb Z}\setminus S)$ and $n+1\in S$, then $f((n,x))=(n+1,x)$ for each $x\in[0,1)$ and $f((n,x))=(n+1,x-1)$ for each $x\in[2,3)$. 
\end{enumerate}

\noindent
By definition we have $f^n(X_m)=X_{m+n}$ for any integers $m$ and $n$. Suppose that $n\in({\mathbb Z}\setminus S)$. Since $X_0=\{0\}\times[0,2)$ is connected and $f^n(X_0)=X_n=\{n\}\times([0,1)\cup[2,3))$ is not connected, we see that $f^n$ is discontinuous. Therefore $n$ is not an element of ${\mathcal C}(X,f)$. Suppose that $n\in S$. For each $m\in({\mathbb Z}\setminus S)$ we see that $f^n$ maps $X_m=\{m\}\times([0,1)\cup[2,3))$ onto $X_{m+n}$. If $m+n\in S$ then $X_{m+n}=\{m+n\}\times[0,2)$ and $f^n((m,x))=(m+n,x)$ for each $x\in[0,1)$ and $f^n((m,x))=(m+n,x-1)$ for each $x\in[2,3)$. Therefore $f^n$ maps $X_m$ continuously onto $X_{m+n}$. 
If $m+n\in({\mathbb Z}\setminus S)$ then $X_{m+n}=\{m+n\}\times([0,1)\cup[2,3))$ and $f^n((m,x))=(m+n,x)$ for each $x\in([0,1)\cup[2,3))$. 
Therefore $f^n$ maps $X_m$ homeomorphically onto $X_{m+n}$. 
Thus we see that $f^n|_{X_m}$ is continuous for each $m\in({\mathbb Z}\setminus S)$. 
Suppose that $m$ is an element of $S$. Then $X_m=\{m\}\times[0,2)$. 
Since $S$ is a submonoid of ${\mathbb Z}$ we see that $m+n$ is also an element of $S$. 
Therefore $X_{m+n}=\{m+n\}\times[0,2)$. We see that $f^n((m,x))=(m+n,x)$ for each $x\in[0,2)$. 
Therefore $f^n$ maps $X_m$ homeomorphically onto $X_{m+n}$. 
Thus we see that $f^n|_{X_m}$ is continuous for each $m\in S$. 
Therefore $f^n$ is continuous. Therefore $n$ is an element of ${\mathcal C}(X,f)$. Thus we have $S={\mathcal C}(X,f)$ as desired. 
$\Box$

\vskip 3mm

\begin{figure}[htbp]
      \begin{center}
\scalebox{0.4}{\includegraphics*{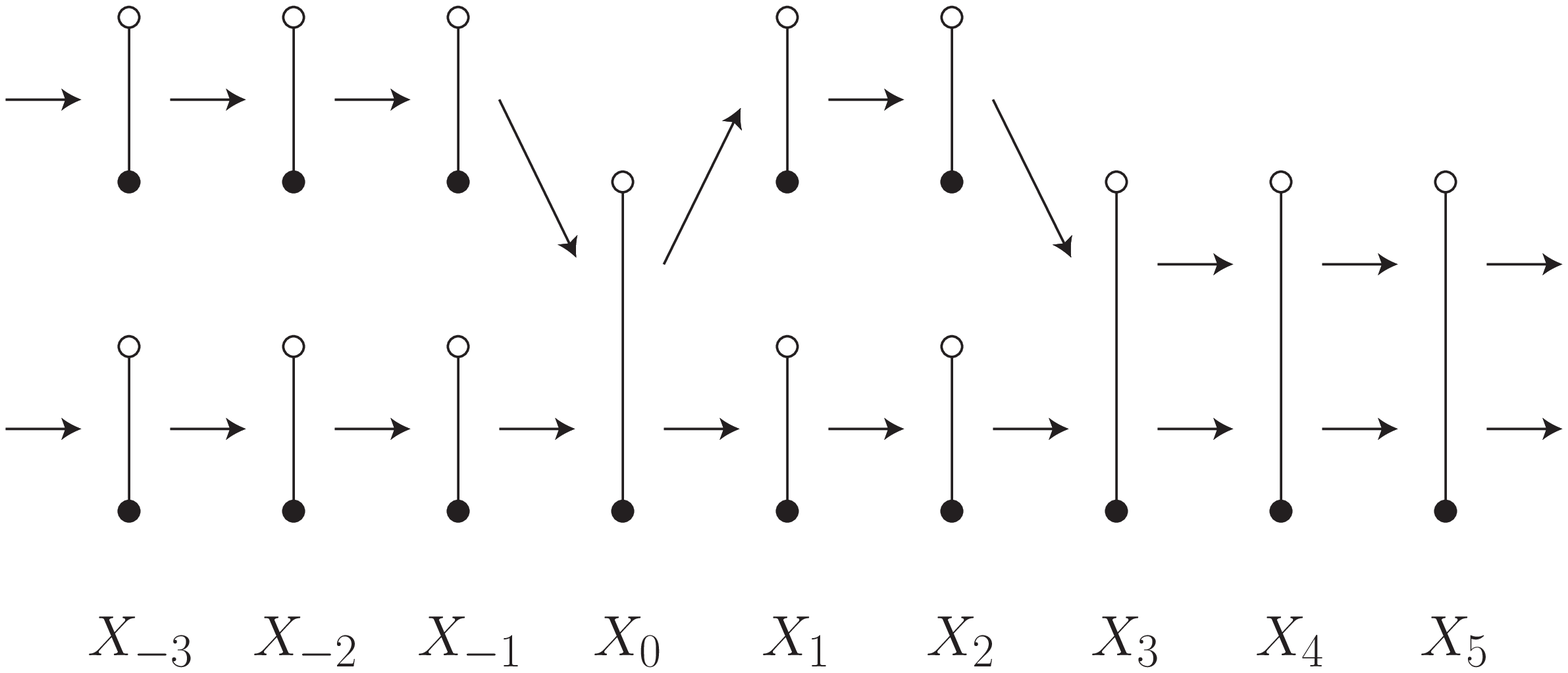}}
      \end{center}
   \caption{}
  \label{bijection}
\end{figure}

\begin{Example}\label{example}
{\rm 
Figure \ref{bijection} illustrates $X$ and $f:X\to X$ in the proof of Theorem \ref{main-theorem} where $S={\mathcal C}(X,f)=\{0\}\cup\{n\in{\mathbb Z}|n\geq3\}$. 
}
\end{Example}

%

We note that the topological type of the topological space $X$ in the proof of Theorem \ref{main-theorem} is independent of the choice of the subset $S$ of ${\mathbb Z}$. Actually $X$ is a disjoint union of countably many semi-open intervals. Thus we have shown the following proposition.

\begin{Proposition}\label{proposition2}
Let $X$ be a disjoint union of countably many semi-open intervals. Then for any submonoid $S$ of ${\mathbb Z}$ there is a bijection $f:X\to X$ such that $S={\mathcal C}(X,f)$. 
\end{Proposition}

We note that not all topological spaces have such a property as $X$ in Proposition \ref{proposition2}. For example, let $X$ be a compact Hausdorff space. 
Then a continuous bijection from $X$ to $X$ is a homeomorphism. 
Therefore, for any bijection $f:X\to X$ the set ${\mathcal C}(X,f)$ is invariant under the map $r:{\mathbb Z}\to{\mathbb Z}$ defined by $r(x)=-x$, $x\in{\mathbb Z}$.

\section{Generalization} 

In this section we reformulate and generalize Theorem \ref{main-theorem} as follows. Let $G$ be a group and $e$ the unit element of $G$. A subset $S$ of $G$ is called a {\it submonoid} of $G$ if $S$ satisfies the following two conditions. 

\begin{enumerate}
\item[(1)]$S$ contains $e$, 
\item[(2)]if $S$ contains $a$ and $b$ then $S$ contains $ab$. 
\end{enumerate}

\noindent
Let $X$ be a topological space. By ${\mathcal B}(X)$ we denote the set of all bijections from $X$ to $X$. Then ${\mathcal B}(X)$ forms a group under the composition of maps. Let $A(X)$ be a subgroup of ${\mathcal B}(X)$. By ${\mathcal C}(A(X))$ we denote the set of all continuous bijections in $A(X)$. Since ${\rm id}_X:X\to X$ is continuous and the composition of two continuous maps is continuous, we see that ${\mathcal C}(A(X))$ is a submonoid of $A(X)$. Let $G$ and $H$ be groups and $S$ and $T$ submonoids of $G$ and $H$ respectively. We say that the pair $(G,S)$ is {\it isomorphic} to the pair $(H,T)$ if there is a group isomorphism $h:G\to H$ such that $h(S)=T$. 

\vskip 3mm

\begin{Theorem}\label{generalization}
Let $G$ be a group and $S$ a submonoid of $G$. Then there is a topological space $X$ and a subgroup $A(X)$ of ${\mathcal B}(X)$ such that the pair $(G,S)$ is isomorphic to the pair $(A(X),{\mathcal C}(A(X)))$. 
\end{Theorem}

\vskip 3mm

\noindent{\bf Proof.}
Let $G$ be a group and $S$ a submonoid of $G$. We give a discrete topology to $G$. Let ${\mathbb R}$ be the $1$-dimensional Euclidean space and $G\times{\mathbb R}$ the product topological space. 
For each element $n$ in $G$ we define a subspace $X_n$ of $G\times{\mathbb R}$ as follows. 
\[
X_n=\left\{ \begin{array}{ll}
\{n\}\times[0,2) & (n\in S)\\
\{n\}\times([0,1)\cup[2,3)) & (n\in (G\setminus S)).\\
\end{array} \right.
\]
Let $\displaystyle{X=\bigcup_{n\in G}X_n}$. Then $X$ is a topological subspace of $G\times{\mathbb R}$. 
For each element $n$ in $G$ we define a bijection $f_n:X\to X$ by the followings. 

\begin{enumerate}
\item[(1)]if $m,mn\in S$, then $f_n((m,x))=(mn,x)$ for each $x\in[0,2)$, 
\item[(2)]if $m,mn\in(G\setminus S)$, then $f_n((m,x))=(mn,x)$ for each $x\in([0,1)\cup[2,3))$, 
\item[(3)]if $m\in S$ and $mn\in(G\setminus S)$, then $f_n((m,x))=(mn,x)$ for each $x\in[0,1)$ and $f_n((m,x))=(mn,x+1)$ for each $x\in[1,2)$, 
\item[(4)]if $m\in(G\setminus S)$ and $mn\in S$, then $f_n((m,x))=(mn,x)$ for each $x\in[0,1)$ and $f_n((m,x))=(mn,x-1)$ for each $x\in[2,3)$. 
\end{enumerate}

\noindent
For any two elements $m$ and $n$ in $G$ we see by definition that $f_n\circ f_m=f_{mn}$. 
Let $A(X)$ be the subgroup of ${\mathcal B}(X)$ defined by $A(X)=\{f_n|n\in G\}$. Then we see that the group $A(X)$ is isomorphic to the group $G$. 
Then by an entirely analogous argument as in the proof of Theorem \ref{main-theorem} we see that ${\mathcal C}(A(X))=\{f_n|n\in S\}$. 
Thus we see that the pair $(A(X),{\mathcal C}(A(X)))$ is isomorphic to the pair $(G,S)$ as desired. 
$\Box$

\vskip 3mm

\begin{Remark}\label{remark}
{\rm
(1) In general the group ${\mathcal B}(X)$ is so big that we should take a subgroup $A(X)$ of ${\mathcal B}(X)$ as in the statement of Theorem \ref{generalization}. 
In fact there is a group $G$ that is not isomorphic to ${\mathcal B}(X)$ for any set $X$. 
For example, it is easy to check that ${\mathcal B}(X)$ is not isomorphic to a cyclic group of order $3$ for any set $X$. 

\noindent
(2) Even in the case that a group $G$ is isomorphic to ${\mathcal B}(X)$ for some set $X$, not all pair $(G,S)$ is realized by the pair $({\mathcal B}(X),{\mathcal C}({\mathcal B}(X)))$ under any topology on $X$. Let $G=S_3$ be a symmetric group of degree $3$. Note that every submonoid of a finite group $G$ is a subgroup of $G$. We will see that the pair $(S_3,C_3)$ is not realized where $C_3$ is a cyclic group of order $3$. 
It is clear that ${\mathcal B}(X)$ is isomorphic to $S_3$ if and only if $X$ contains exactly $3$ points. Therefore we may suppose without loss of generality that $X=\{a,b,c\}$. Then, up to self-homeomorphism, there are $9$ topologies on $X$. They are ${\mathcal D}_1=\{\emptyset,X\}$, ${\mathcal D}_2=\{\emptyset,\{a\},X\}$, ${\mathcal D}_3=\{\emptyset,\{a,b\},X\}$, ${\mathcal D}_4=\{\emptyset,\{a\},\{a,b\},X\}$, ${\mathcal D}_5=\{\emptyset,\{a\},\{b,c\},X\}$, ${\mathcal D}_6=\{\emptyset,\{a\},\{b\},\{a,b\},X\}$, ${\mathcal D}_7=\{\emptyset,\{a\},\{a,b\},\{a,c\},X\}$, ${\mathcal D}_8=\{\emptyset,\{a\},\{b\},\{a,b\},\{a,c\},X\}$ and ${\mathcal D}_9=2^X$. Then we see that the subgroup ${\mathcal C}({\mathcal B}(X,{\mathcal D}_i))$ of ${\mathcal B}(X,{\mathcal D}_i)$ is the trivial group for $i=4,8$, a cyclic group of order $2$ for $i=2,3,5,6,7$ and the symmetric group of degree $3$ ${\mathcal B}(X,{\mathcal D}_i)$ for $i=1,9$. Thus ${\mathcal C}({\mathcal B}(X,{\mathcal D}_i))$ is not a cyclic group of order $3$ for any $i$. 
}
\end{Remark}

\vskip 3mm

Next we give a variation of Theorem \ref{generalization} as follows. A {\it monoid} $M$ is a semigroup with the unit element $e$. Namely $M$ has an associative binary operation such that $xe=ex=x$ for any element $x\in M$. A subset $S$ of a monoid $M$ is said to be a {\it submonoid} of $M$ if $e$ is an element of $S$ and for any elements $a$ and $b$ of $S$ the element $ab$ is an element of $S$. Let $X$ be a topological space. By ${\mathcal M}(X)$ we denote the set of all maps from $X$ to $X$. Then ${\mathcal M}(X)$ forms a monoid under the composition of maps. Let $A(X)$ be a submonoid of ${\mathcal M}(X)$. By ${\mathcal C}(A(X))$ we denote the set of all continuous maps in $A(X)$. Then we see as before that ${\mathcal C}(A(X))$ is a submonoid of $A(X)$. Let $M$ and $N$ be monoids and $S$ and $T$ submonoids of $M$ and $N$ respectively. We say that the pair $(M,S)$ is {\it isomorphic} to the pair $(N,T)$ if there is a monoid isomorphism $h:M\to N$ such that $h(S)=T$. 

\vskip 3mm

\begin{Theorem}\label{generalization2}
Let $M$ be a monoid and $S$ a submonoid of $M$. Then there is a topological space $X$ and a submonoid $A(X)$ of ${\mathcal M}(X)$ such that the pair $(M,S)$ is isomorphic to the pair $(A(X),{\mathcal C}(A(X)))$. 
\end{Theorem}

\vskip 3mm

\noindent{\bf Proof.} 
We define a topological space $X$ to be a subspace of $M\times{\mathbb R}$ as in the proof of Theorem \ref{generalization}. The map $f_n:X\to X$ is also defined in the same way for each element $n$ of $M$. The only difference is that the map $f_n$ is not a bijection in general. Note that in the proof of Theorem \ref{generalization} the assumption that $n$ has an inverse element in the group $G$ assured the fact that $f_n:X\to X$ is a bijection. Then the rest of the proof is the same as the proof of Theorem \ref{generalization}. 
$\Box$

\vskip 3mm

Finally we give another variation of Theorem \ref{generalization} as follows. As we have already remarked, if $X$ is compact Hausdorff and $f:X\to X$ is a continuous bijection, then $f^{-1}:X\to X$ is also continuous. Therefore for any subgroup $A(X)$ of ${\mathcal B}(X)$ the submonoid ${\mathcal C}(A(X))$ of $A(X)$ is a subgroup of $A(X)$. Then we have the following theorem. 

\vskip 3mm

\begin{Theorem}\label{generalization3}
Let $G$ be a group and $H$ a subgroup of $G$. Then there is a compact Hausdorff space $X$ and a subgroup $A(X)$ of ${\mathcal B}(X)$ such that the pair $(G,H)$ is isomorphic to the pair $(A(X),{\mathcal C}(A(X)))$. 
\end{Theorem}

\vskip 3mm

\noindent{\bf Proof.} 
We give a discrete topology to $G$. Let $G\times[0,1]$ be the product topological space and $X=(G\times[0,1])\cup\{\infty\}$ the one-point compactification of $G\times[0,1]$. Then $X$ is a compact Hausdorff space. For each element $n$ in $G$ we define a bijection $f_n:X\to X$ by the followings. 

\begin{enumerate}
\item[(1)]For each $m$ in $G$ and $x$ in $(0,1)$, $f_n((m,x))=(mn,x)$. 
\item[(2)]If $m,mn\in H$ or $m,mn\in (G\setminus H)$, then $f_n((m,0))=(mn,0)$ and $f_n((m,1))=(mn,1)$. 
\item[(3)]If $m\in H$ and $mn\in(G\setminus H)$, or $m\in (G\setminus H)$ and $mn\in H$, then $f_n((m,0))=(mn,1)$ and $f_n((m,1))=(mn,0)$. 
\item[(4)]$f_n(\infty)=\infty$. 
\end{enumerate}

\noindent
We see by the definition that the composition $f_n\circ f_m$ is equal to $f_{mn}$ for any elements $m$ and $n$ in $G$. 
Let $A(X)$ be the subgroup of ${\mathcal B}(X)$ defined by $A(X)=\{f_n|n\in G\}$. Then we see that the group $A(X)$ is isomorphic to the group $G$. 
We will show that ${\mathcal C}(A(X))=\{f_n|n\in H\}$. 
First we will show that $f_n$ is continuous at $\infty$ for any $n$ in $G$. Let $U$ be an open neighbourhood of $\infty$. Then $X\setminus U$ is a compact subset of $G\times [0,1]$. Therefore there is a finite subset $F$ of $G$ such that $X\setminus U$ is contained in $F\times[0,1]$. 
Let $V=X\setminus((Fn^{-1})\times[0,1])$. Then $V$ is an open neighbourhood of $\infty$ such that $f_n(V)=X\setminus(F\times[0,1])$ is contained in $U$ as desired. 
Suppose that $n\in(G\setminus H)$. Then $f_n$ maps $\{e\}\times[0,1]$ to $\{n\}\times[0,1]$. Since the unit element $e$ is in $H$, $f_n((e,0))=(n,1)$ and $f_n((e,1))=(n,0)$. Therefore $f_n$ is not continuous and $f_n$ is not in ${\mathcal C}(A(X))$. Suppose that $n\in H$. Let $m$ be an element of $G$. Then we see that $mn$ is an element of $H$ if and only if $m$ is an element of $H$. Therefore $f_n$ maps $\{m\}\times[0,1]$ to $\{mn\}\times[0,1]$ by the formula $f_n((m,x))=(mn,x)$ for each $x$ in $[0,1]$. Therefore the restriction map $f_n|_{\{m\}\times[0,1]}$ is continuous for each $m$ in $G$. Therefore $f_n$ is an element of ${\mathcal C}(A(X))$. Thus the pair $(G,H)$ is isomorphic to the pair $(A(X),{\mathcal C}(A(X)))$. 
$\Box$

\vskip 3mm

\begin{Remark}\label{remark}
{\rm
Theorem \ref{generalization}, Theorem \ref{generalization2} and Theorem \ref{generalization3} concern the pairs $(G,S)$, $(M,S)$ and $(G,H)$ respectively. 
There are some known results not on a pair but on a single group or a single monoid. 
It is shown in \cite{Groot} that for any group $H$ there exists a topological space $X$ such that the group of all self-homeomorphisms of $X$ is isomorphic to $H$. 
It is shown in \cite{Trnkova} that for any monoid $S$ there exists a topological space $X$ such that the monoid of all nonconstant continuous maps from $X$ to $X$ is isomorphic to $S$. 
}
\end{Remark}

\vskip 3mm

\vskip 3mm

\section*{Acknowledgments} The author is grateful to Professor Kazuhiro Kawamura for his helpful comments.

\vskip 3mm

\vskip 3mm

{
 
}

\vskip 3mm

\vskip 3mm

\end{document}